\documentclass[twoside]{amsart}
\usepackage{latexsym}
\usepackage{amssymb,amsmath,amsopn}
\usepackage[dvips]{graphicx}   %To insert ps figures
\usepackage{color,epsfig}      %To insert ps figures 

\newtheorem{thm}{Theorem}

\theoremstyle{definition}

\renewcommand{\Re}{\mathbb R}

\renewcommand{\S}{\mathbb{S}}

\renewcommand{\P}{\mathcal{P}}

\DeclareMathOperator{\ind}{Ind}

\DeclareMathOperator{\In}{In}

\begin{document}
\title[On universal covers]{On universal covers for four-dimensional sets of a given diameter}

\author[Z. L\'angi]{Zsolt L\'angi}
                                                                               
\address{Zsolt L\'angi, Dept. of Geometry, Budapest University of Technology,
Budapest, Egry J\'ozsef u. 1., Hungary, 1111}
\email{zlangi@math.bme.hu}
                                                                        

\begin{abstract}
Makeev proved that among centrally symmetric four-dimensional polytopes, with more than twenty facets and circumscribed about the Euclidean ball of diameter one, there is no universal cover for the family of unit diameter sets.
In this paper we examine the converse  problem, and prove that each centrally symmetric polytope, with at most fourteen facets and circumscribed about the Euclidean ball of diameter one, is a universal cover for the family of unit diameter sets.
\end{abstract}
\maketitle

A convex body $C$ in the Euclidean $n$-space $\Re^n$ is called a
\emph{universal cover for sets of diameter $d$},
if for every set of diameter $d$ there is a congruent copy of $C$ containing it.
The problem of finding universal covers for sets of a given diameter, or equivalently, for unit diameter sets,
is a long-standing question of discrete geometry.
These universal covers are used, in particular, for the solution of Borsuk's problem, that asks the minimal number of
subsets of smaller diameters that an $n$-dimensional set can be partitioned into.
For information about Borsuk's problem and its relationship with universal covers, the reader is referred to \cite{R07}.

As a special case, we may consider universal covers in the family $\P_n$ of
centrally symmetric $n$-dimensional polytopes circumscribed about a Euclidean ball of diameter one.
Since every set in the Euclidean space is contained in a constant-width body of the same diameter, we may rephrase
this problem as finding polytopes that can be circumscribed about any $n$-dimensional body of constant width one.
This problem is related to Knaster's problem, that asks which finite point sets
on the Euclidean sphere $\S^m$ has the property that for any continuous function
$f : \S^m \to \Re^n$, $f$ is constant on a congruent copy of $S$ (cf., e.g. \cite{K99}).

Let $D_n$ denote the dual of the difference body of an $n$-dimensional regular simplex such that
$D_n$ is circumscribed about the Euclidean ball of unit diameter.
In 1994, Makeev \cite{M94} conjectured that $D_n$ is a universal cover for unit diameter sets.
This conjecture is partly motivated by the fact that $D_n$ has $n(n+1)$ facets, and by a result of Makeev
\cite{M94} that no universal cover in $\P_n$ has more than $n(n+1)$ facets.
The $n=2$ case of Makeev's conjecture is known as P\'al's lemma, and has been known since 1920 (cf. \cite{P20}).
The $n=3$ case was proven independently by Makeev \cite{M97}, by G. Kuperberg \cite{K99} and by Hausel, Makai, Jr. and Sz\H ucs \cite{HMSz00}.
In \cite{M97} and \cite{K99}, the main idea of the proofs is that, homologically,
for any convex body of constant width one, there are an odd number of congruent copies of $D_3$
circumscribing it, which yields that geometrically there is at least one.
The author of \cite{K99} remarks also that in dimension four, homologically, there is zero circumscribed copy of not only $D_4$, which has twenty facets, but also of the regular cross-polytope, which has sixteen facets.

Our main result is the following. We note that this result, in some sense, is converse to the result in \cite{M94}, mentioned in the previous paragraph.

\begin{thm}\label{thm:main}
Every polytope $P \in \P_4$, with at most fourteen facets, is a universal cover for unit diameter sets.
\end{thm}

To prove this theorem, first we introduce a topological invariant, called Smith index,
and recall some estimates regarding it.

Let $T$ be a fixed-point-free involution, with period two, defined on the topological spaces $X$ and $Y$
(for the terminology, cf., for example \cite{ES52}).
Then the pairs $(X,T)$ and $(Y,T)$ are called \emph{$T$-spaces},
and a continuous mapping $f : X \to Y$ with $Tf = fT$ is called a \emph{$T$-map}.

Let $S(X)$ be the singular chain complex of $X$, with $Z_2$ as the set of coefficients, and
let $\rho=I-T=I+T$, where $I$ is the identity operator on $S(X)$.
Then, as $T$, and thus also $\rho$, commutes with the boundary operator $\partial$ of $S(X)$, we have that $\rho$ is a chain homomorphism of $S(X)$ into itself.
We denote the image of this homomorhism by $S^\rho(X)$.
Composing $S^\rho(X)$ with the homology theory functor of $Z_2$, we obtain the singular Smith homology groups with coefficient group $Z_2$.
We denote the $k$th group of this theory by $H_k^\rho(X)$, and the corresponding homology group of $S(X)$ by $H_k(X)$.
Theorem 1.6 of \cite{G61} states that if $(X,T)$ is a Haussdorff $T$-space and $X^*$ is its orbit space, then $H_k^\rho(X,T) \approx H_k(X^*,Z_2)$,
for every value of $k$, where $H_k(X^*,Z_2)$ is the ordinary $k$th singular homology group of $X^*$.

Now, let $i$ denote the inclusion operator, and $\Delta_k$ be the boundary operator of $H_k^\rho(X)$.
Then
\begin{equation}\label{eq:exact}
\ldots \xleftarrow{i} H_{k-1}^\rho(X) \xleftarrow{\Delta_k} H_k^\rho(X) \xleftarrow{\rho} H_k(X) \xleftarrow{i} H_k^\rho(X) \xleftarrow{\Delta_{k+1}} \ldots
\end{equation}
is an exact sequence.
This follows from the fact that $Z_2$ is a field, and thus, the couple $(i,\rho)$ is a direct couple, and hence the observation follows from Theorem 2.7, p. 128 of \cite{ES52}.

Furthermore, we may define a homomorphism
\[
\tilde{\Delta}_0 : H_0^\rho(X) \to Z_2 \hbox{ by}
\]
\[
\tilde{\Delta}_0 = \In \circ \rho^{-1},
\]
where $\In$ is the Kronecker index homomorphism of the ordinary singular homology theory.
This is well defined, as the kernel of $\rho$ is contained in the kernel of $\In$.
Then, we define
\[
\ind : H_n^\rho(X) \to Z_2 \hbox{ by}
\]
\[
\ind = \tilde{\Delta}_0 \circ \Delta_1 \circ \Delta_2 \circ \ldots \circ \Delta_n.
\]

Finally, the Smith index $\ind(\emptyset,T)$ is zero. If $X$ is not empty, then $\ind(X,T)$ is the largest integer
such that $\ind(H_n^\rho(X)) \neq 0$, if it exists, and otherwise it is defined to be infinity.

An important property of Smith index is stated in the following theorem, proved by Geraghty
(cf. the remark after Lemma 2.2 of \cite{G61}).

\begin{thm}[Geraghty, 1961]\label{thm:Smithimp}
If there is a $T$-map $f : (X,T) \to  (Y,T)$, then
$\ind(X,T) \leq \ind(Y,T)$.
\end{thm}%

Another important observation that we use later is that
if $H_0(X) = Z_2$, and $H_i(X) = 0$ for $i=2,3,\ldots,n-1$, then $\ind(X,T)\geq n$.
Indeed, since $H_0(X) = Z_2$ and since the sequence in (\ref{eq:exact}) is exact, we have $H_0^\rho(X)=Z_2$.
Thus $\Delta_1$ is onto, and as $H_i(X) = 0$ for $i=2,3,\ldots,n-1$, we obtain that $\Delta_i$ is onto for $i=2,3,\ldots,n$.
Since a zero-dimensional $\rho$-cycle consisting of a point and its $T$-image has index one, it yields that $\ind(X,T) \geq n$.
In particular, it is well-known that for the Euclidean sphere $\S^n$ with the usual antipodal mapping $T$, $\ind (\S^n, T) = n$.

Now recall the notion of Stiefel manifold; that is,
the topological space of the orthonormal $k$-frames in $\Re^n$, denoted by $V_{n,k}$.
Observe that there is a natural homeomorphism between $V_{n,n-1}$ and $SO(n)$, and between $V_{n,1}$ and $\S^{n-1}$.
We use the following estimates, proved by Geraghty, regarding the Smith indices of Stiefel manifolds
with respect to the ususal antipodal mapping $T$ (cf. \cite{G61}).

\begin{thm}[Geraghty, 1961]\label{thm:estimate}
If $s$ is the largest power of $2$ that divides $2n$, then 
\[
s-1 \leq \ind(SO(2n),T) \leq \ind(\S^{2n-1},T) = 2n-1.
\]
In particular, if $n$ is a power of $2$, then $\ind(SO(2n),T) = 2n-1$.
\end{thm}

For completeness, we recall the proof of this estimate from \cite{G61}.

\begin{proof}
Clearly, it is sufficient to prove the general estimate.
Consider the Stiefel manifolds $V_{2n,k}$ with the usual antipodal mapping as $T$.
Observe that by deleting the $k$th member of the frame, we have a sequence of $T$-maps
\[
(SO(2n),T) = (V_{2n,2n-1},T) \to (V_{2n,2n-2},T) \to \ldots, \to (V_{2n,1},T) = (\S^{2n-1},T).
\]
The orbit space of $SO(2n)$ is the projective special orthogonal group $PSO(2n)$. Thus, the singular Smith homology
groups of $SO(2n)$ are the ordinary singular homology groups of $PSO(2n)$.

The Poincar\'e polynomial of $PSO(2n)$, with the coefficients in $Z_2$, is
\[
P(t)=(1+t)(1+t^2)\ldots(1=t^{s-1}) \cdot (1+t+\ldots+t^s+t^s)(1+t^{s+1})\ldots(1+t^{2n-1}),
\]
where $s$ is the largest power of $2$ dividing $2n$ (cf. \cite{B54}).

Furthermore, the Poincar\'e polynomial of $SO(2n)$ (cf. \cite{B51}), with the coefficients in $Z_2$, is
\[
Q(t)=(1+t)(1+t^2)\ldots(1+t^{2n-1}).
\]

Now, let the coefficient of $t^i$ in P(t), or in other words the Betti number of $H_i(PSO(2n),Z_2)=H_i^\rho(SO(2n),T)$,
be denoted by $B_i^\rho$, and similarly, let $B_i$ denote the coefficient of $t^i$ in $Q(t)$.
Then for $i=0,1,\ldots,s-1$, we have
\[
B_i^\rho = \sum_{j=0}^i B_j.
\]

Thus, in terms of the Betti numbers, the sequence in (\ref{eq:exact}) is
\[
0 \leftarrow B_0 \xleftarrow{\rho} B_o \xleftarrow{\Delta_1} B_0+B_1 \xleftarrow{\rho} B_1 \xleftarrow{i} B_0+B_1 \xleftarrow{\Delta_2} B_0 + B_1 + B_2 \xleftarrow{\rho} \ldots
\]
and hence, $\Delta_i$ is onto for $i=1,2,\ldots,s-1$. Since $H_0^\rho(SO(2n),T) = Z_2$, we have $\ind(SO(2n),T) \geq s-1$.
Since $\ind(\S^{2n-1},T)=2n-1$, and since $T$-maps do not decrease the value of Smith index, we have that
\[
s-1 \leq \ind(SO(2n),T) \leq \ind(\S^{2n-1},T) = 2n-1.
\]
\end{proof}

Now we are ready to prove Theorem~\ref{thm:main}.

\begin{proof}
Consider a convex body $C \subset \Re^4$ of constant width one.
We may assume that $P$ has exactly fourteen facets.
Let $K_1, K_2, \ldots, K_7$ denote the seven infinite strips bounded by pairs of parallel
facet-hyperplanes of $P$, and note that $P = \bigcap_{i=1}^{7} K_i$.
Observe that the width of any of these strips is one.
Note also that in any system of vectors spanning $\Re^4$, there are four that also span $\Re^4$.
Applying this observation for the normal vectors of the facet-hyperplanes of $P$, we obtain
that, among $K_1, K_2, \ldots, K_7$, there are four strips such that the hyperplanes bisecting them
intersect in a singleton.
We may assume that these strips are $K_1, K_2, K_3$ and $K_4$, which yields the (unique) existence of a
translation vector $x$ such that $C \subset x+ \bigcap_{i=1}^4 K_i$.

Consider an arbitrary element $\tau \in SO(4)$.
We define a function $g : SO(4) \to \Re^3$ in the following way.
Let $x_\tau$ denote the unique vector with the property that $C \subset x_\tau+ \bigcap_{i=1}^4 \tau(K_i)$.
Then the three coordinates of $g(\tau)$ are the signed distances, from $x_\tau$, of the three hyperplanes
bisecting $\tau(K_5)$, $\tau(K_6)$ and $\tau(K_7)$.

Let $T$ denote the usual antipodal mapping, and note that by Theorem~\ref{thm:estimate}, we have $\ind(SO(4),T)=3$.
As by Theorem~\ref{thm:Smithimp}, $T$-maps do not decrease the value of Smith index, and as $\ind(\S^2,T)=2$, there is no $T$-map from $(SO(4),T)$ to $(\S^2,T)$.
Similarly like in the proof of the classical Borsuk-Ulam theorem, from this it can be shown that for any $T$-map $g:SO(4) \to \Re^3$, there is a point mapped to the origin
(this property follows also from Theorem 4.2 of \cite{G61}).
\end{proof}

We note that in our consideration, we have shown also the following, more general statement.

\begin{thm}
Every polytope $P \in \P_{2m}$ with at most $2m+2\ind(SO(2m),T)$ facets is a universal cover for
sets of diameter one.
\end{thm}

{\bf Acknowledgements.}
The author is indebted to M. Nasz\'odi and an anonymous referee for their helpful remarks.

\end{document}